\documentclass[11pt]{article}
\usepackage{latexsym, amsfonts}
\newcommand{\be}{\begin{equation}}
\newcommand{\ee}{\end{equation}}
\newcommand{\bea}{\begin{eqnarray}}
\newcommand{\eea}{\end{eqnarray}}
\newcommand{\bean}{\begin{eqnarray*}}
\newcommand{\eean}{\end{eqnarray*}}
\newcommand{\brray}{\begin{array}}
\newcommand{\erray}{\end{array}}
\newcommand{\ben}{\begin{equation}{nonumber}}
\newcommand{\een}{\end{equation}{nonumber}}

\newtheorem{dfn}{Definition}[section]
\newtheorem{thm}[dfn]{Theorem}
\newtheorem{lmma}[dfn]{Lemma}
\newtheorem{ppsn}[dfn]{Proposition}
\newtheorem{crlre}[dfn]{Corollary}
\newtheorem{xmpl}[dfn]{Example}
\newtheorem{rmrk}[dfn]{Remark}
\newcommand{\bdfn}{\begin{dfn}}
\newcommand{\bthm}{\begin{thm}}
\newcommand{\blmma}{\begin{lmma}}
\newcommand{\bppsn}{\begin{ppsn}}
\newcommand{\bcrlre}{\begin{crlre}}
\newcommand{\bxmpl}{\begin{xmpl}}
\newcommand{\brmrk}{\begin{rmrk}}
\newcommand{\edfn}{\end{dfn}}
\newcommand{\ethm}{\end{thm}}
\newcommand{\elmma}{\end{lmma}}
\newcommand{\eppsn}{\end{ppsn}}
\newcommand{\ecrlre}{\end{crlre}}
\newcommand{\exmpl}{\end{xmpl}}
\newcommand{\ermrk}{\end{rmrk}}

\newcommand{\IC}{\mathbb{C}}

\newcommand{\cla}{{\cal A}}
\newcommand{\clb}{{\cal B}}

\newcommand{\cld}{{\cal D}}
\newcommand{\cle}{{\cal E}}

\newcommand{\clg}{{\cal G}}
\newcommand{\clh}{{\cal H}}

\newcommand{\clk}{{\cal K}}

\newcommand{\clm}{{\cal M}}
\newcommand{\cln}{{\cal N}}

\newcommand{\cls}{{\cal S}}

\newcommand{\clv}{{\cal V}}

\def\a*{{\cal A}_{h,*}}
\def\B{{\cal B}(h)}
\def\B1{{\cal B}_1(h)}
\def\b{{\cal B}^{\rm s.a.}(h)}
\def\b1{{\cal B}^{\rm s.a.}_1(h)}

\newcommand{\ot}{\otimes}

\newcommand{\raro}{\rightarrow}

\def \qed {$\Box$}

\begin{document}
\begin{center}
{\large {\bf Some remarks on the action of  quantum isometry groups}}\\
by\\
{\large Debashish Goswami {\footnote {The author gratefully acknowledges support obtained  from the Indian National Academy of Sciences through the grants for a project on `Noncommutative Geometry and Quantum Groups'.}}}\\
{\large Stat-Math Unit, Kolkata Centre,}\\
{\large Indian Statistical Institute}\\
{\large 203, B. T. Road, Kolkata 700 108, India}\\
\end{center}
\begin{abstract}
We give a new  sufficient condition on a spectral triple  to ensure that the quantum group of orientation and volume preserving isometries defined in \cite{qorient} has a $C^*$-action on the underlying $C^*$ algebra. 
\end{abstract}
\section{Introduction}
Taking motivation from the work of Wang, Banica, Bichon and others (see 
\cite{free}, \cite{wang}, 
\cite{ban1}, \cite{ban2}, \cite{bichon}, \cite{univ1} and references therein), we have given a definition 
 of quantum isometry group based on a `Laplacian' in \cite{goswami}, and then followed it up by a formulation of `quantum group of orientation preserving isometries' in \cite{qorient} (see also \cite{qdisc}, \cite{q_sphere}, \cite{jyotish}, \cite{jyotish_1} for many explicit computations). 
The main result of \cite{qorient} is that given a spectral triple (of compact type) $(\cla^\infty, \clh, D)$  and a positive unbounded operator $R$
 commuting with $D$,  there is a universal object in the category of compact quantum groups which have a unitary representation (say $U$) on the Hilbert space $\clh$ w.r.t. which $D$ is equivariant and the normal (co)-action $\alpha_U$ of the quantum group obtained by conjugation by the unitary representation leaves the weak closure of $\cla^\infty$ invariant and preserves a canonical functional $\tau_R$ called `$R$-twisted volume form' described in \cite{qorient}. The  Woronowicz subalgebra of this  universal quantum group genearted by the `matrix elements' of $\alpha_U(a)$, $a \in \cla^\infty$ is called `the quantum group of orientation and $R$-twisted volume preserving isometries' and  denoted by $QISO^+_R(D)$. 

However, it is not clear from the definition and construction of this quantum group whether $\alpha_U$ is a $C^*$ -action of $QISO^+_R(D)$ on 
the $C^*$ algebra generated by $\cla^\infty$ (in the sense of Woronowicz and Podles). The problem is that to prove  the existence of a universal object in \cite{qorient}
  we had to make use of the Hilbert space and the strong operator topology coming from it. This is the reason why we demanded only the stability of the von Neumann
 algebra generated vy $\cla^\infty$ in the definition of an isometric and orientation preserving quantum group action on the spectral triple 
$(\cla^\infty, \clh, D)$. In a sense, we worked in a suitable category of `measurable' actions and could prove the existence of a universal object there. 
 In the classical situation, i.e. for isometric orientation preserving group actions on Riemannian manifolds, the aparently weaker condition of measurability 
 turns out to be equivalent to a topological, in fact smooth  action, thanks to the Sobolev's theorem.The analogues question in the noncommutative situation is to
 see whether $\alpha_U$ (which is a-prori only a normal action) is a $C^*$ action, and we shall show that under some assumptions which very much resemble the 
  clasical Sobolev's theorem, we can indeed answer this question in the affirmative. In fact, we have already given a number of sufficient conditions 
for ensuring $C^*$-action in \cite{qorient}, and the conditions  given in the present article add to this list, strengthening our belief that the action of quantum
 group of orientation preserving isometries is in general a $C^*$ action.

\section{Preliminaries}
\subsection{Generalities on quantum groups and their action}
We review some basic facts about quantum groups (see, e.g. \cite{woro}, \cite{woro1}, \cite{vandaelenotes} and references therein).  A 
compact quantum group (to be abbreviated as CQG from now on)  is given by a pair $(\cls, \Delta)$, where $\cls$ is a unital separable $C^*$ algebra 
equipped
 with a unital $C^*$-homomorphism $\Delta : \cls \raro \cls \otimes \cls$ (where $\otimes$ denotes the injective tensor product)
  satisfying \\
  (ai) $(\Delta \ot id) \circ \Delta=(id \ot \Delta) \circ \Delta$ (co-associativity), and \\
  (aii) the linear span of $ \Delta(\cls)(\cls \ot 1)$ and $\Delta(\cls)(1 \ot \cls)$ are norm-dense in $\cls \ot \cls$. \\
  It is well-known (see \cite{woro}, \cite{woro1}) that there is a canonical dense $\ast$-subalgebra $\cls_0$ of $\cls$, consisting of the matrix coefficients of
   the finite dimensional unitary (co)-representations (to be defined shortly) of $\cls$, and maps $\epsilon : \cls_0 \raro \IC$ (co-unit) and
   $\kappa : \cls_0 \raro \cls_0$ (antipode)  defined
    on $\cls_0$ which make $\cls_0$ a Hopf $\ast$-algebra.

    We say that  the compact quantum group $(\cls,\Delta)$ (co)-acts on a unital $C^*$ algebra $\clb$,
    if there is a  unital $C^*$-homomorphism (called an action) $\alpha : \clb \raro \clb \ot \cls$ satisfying the following :\\
    (bi) $(\alpha \ot id) \circ \alpha=(id \ot \Delta) \circ \alpha$, and \\
    (bii) the linear span of $\alpha(\clb)(1 \ot \cls)$ is norm-dense in $\clb \ot \cls$.\\ 
 It is known that the above is equivalent to the existence of a dense unital $\ast$-subalgebra $\clb_0$ of $\clb$ on which the action $\alpha$ is an algebraic action of the Hopf algebra $\cls_0$, i.e. $\alpha$ maps $\clb_0$ into $\clb_0 \ot_{\rm alg} \cls_0$ and also $({\rm id} \ot \epsilon) \circ \alpha={\rm id}$ on $\clb_0$.

Such an action will be called a $C^*$ or topological  action, to distinguish it 
 from a nornal action on a von Neumann algebra, which we briefly mention later.
     \vspace{1mm}\\

  \bdfn
  
   A unitary ( co ) representation of a compact quantum group $ ( S, \Delta ) $ on a Hilbert space $ \clh $ is a map $ U $ from $ \clh $ to the Hilbert $\cls$ module $ \clh \otimes \cls $  such that the  element $ \widetilde{U} \in \clm ( \clk ( \clh ) \otimes \cls ) $ given by $\widetilde{U}( \xi \ot b)=U(\xi)(1 \ot b)$ ($\xi \in \clh, b \in \cls)$) is a unitary satisfying  $$ ({\rm  id} \otimes \Delta ) \widetilde{U} = {\widetilde{U}}_{(12)} {\widetilde{U}}_{(13)},$$ where for an operator $X \in \clb(\clh_1 \ot \clh_2)$ we have denoted by $X_{12}$ and $X_{13}$ the operators $X \ot I_{\clh_2} \in \clb(\clh_1 \ot \clh_2 \ot \clh_2)$, and $\Sigma_{23} X_{12} \Sigma_{23}$ respectively ($\Sigma_{23}$ being the unitary on $\clh_1 \ot \clh_2 \ot \clh_2$ which flips the two copies of $\clh_2$).
  
Given a unitary representation $U$ we shall denote by $\alpha_U$ the $\ast$-homomorphism $\alpha_U(X)=\widetilde{U}(X \ot 1){\widetilde{U}}^*$ for $X \in \clb(\clh)$. For a  not necessarily bounded, densely defined (in the weak operator topology)  linear functional $\tau$ on $\clb(\clh)$,  we say that $\alpha_U$ preserves $\tau$ if $\alpha_U$ maps a suitable (weakly) dense $\ast$-subalgebra   (say $\cld$) in the domain of $\tau$ into $\cld \ot_{\rm alg} \cls$ and $( \tau \ot {\rm id}) (\alpha_U(a))=\tau(a)1_\cls$  for all $a \in \cld$. When $\tau$ is bounded and normal, this is equivalent to $(\tau \ot {\rm id}) (\alpha_U(a))=\tau(a) 1_\cls$ for all $a \in \clb(\clh)$. 

We say that a (possibly unbounded) operator $T$ on $\clh$ commutes with $U$ if $T \ot I$ (with the natural domain) commutes with $\widetilde{U}$. Sometimes such an operator will be called $U$-equivariant.
\edfn

 We also need to consider  Hopf von Neumann algebra  or von Neumann algebraic  quantum group of compact type. This is given by a von Neumann algebra
 $\clm$ equipped with a normal coassociative coproduct and also a faithful normal state (Haar state) invariant under the coproduct. We refer to \cite{vnqgp}, \cite{vnqgp2}
 for more details, in fact for a locally compact von Neuann algebaric quantum group, of which those of compact type form a very special and relatively simple 
 class. We shall actually be concerned with the canonical Hopf von Neumann algebra coming from the GNS representation of a CQG (w.r.t. the Haar state). Let us 
 mention that the Haar state (say $h$) on a CQG $\cls$ is not necessarily faithful, but it is faithful on the dense $\ast$-algebra $\cls_0$ mentioned before. 
 Let $\rho_h : \cls \raro \clb(L^2(h))$  be the GNS representation. It is easily seen that the coproduct on $\cls$ is implemented by the canonical left regular unitary
 representation, say $U$, on this space, and then $\alpha_U$ can be used as the definition of a normal coassociative coproduct on the von Neumann algebra 
 $\clm=\rho_h(\cls)^{\prime \prime} \subseteq \clb(L^2(h))$. Since the Haar state is a vector state in $L^2(h)$ and hence normal invariant state, it is clear that 
 $\clm$ is a compact type von Neumann algebraic quantum group. 

 We remark here that the definition of unitary representation as well as 
(co)-action has a natural analogue for such von Neumann algebraic quantum groups, and it can be shown that any  such reprsentation
 decomposes into direct sum of irreducible ones, and that any irreducibe reprsentation of a compact Hopf von Neumann algebra
 is finte dimensional. The proofs of these facts are almost the same as the proof in the $C^*$ case, with the only dfferenece being that the norm topology 
 must be replaced by appropriate  
  weak or strong operator topology. Moreover, we remark that given a $C^*$-action $\alpha : \cla \raro \cla \ot \cls$ of the CQG $\cls$ on a separable unital $C^*$ algebra imbedded in $\clb(\clh)$ for some Hilbert space $\clh$, such that the action $\alpha$ is implemented by some unitary representation $U$ of the CQG $\cls$ on $\clh$, i.e. $\alpha=\alpha_U$,  we can canonically construct a normal action of the Hopf von Neumann algebra $\clm$ as follows. First, replace the action $\alpha$ by $\alpha_h=({\rm id} \ot \rho_h) \circ \alpha$, which is an action of $\rho_h(\cls)$, and note that $\alpha_h(a)=U_h (a \ot 1)U_h^*$ for $a \in \cla$, where $U_h:=({\rm id} \ot \rho_h)(U)$, a unitary representation of $\rho_h(\cls)$ on $\clh$. Now we use the right hand side of the above to extend the definition of $\alpha_h$ to the whole of $\clb(\clh)$, in particular on $\cln:=\cla^{\prime \prime} \subseteq \clb(\clh)$, and observe that this indeed gives a normal action of  $\clm=\rho_h(\cla)^{\prime \prime}$ on $\cln$. 

\subsection{Quantum group of orientation and volume preserving isometries}
Next we give an overview of the definition of quantum isometry groups as in \cite{qorient}.
\bdfn 
         \label{def_q_fam}A quantum family of orientation preserving  isometries for the spectral triple $({\cla^\infty}, \clh, D)$ is given by a pair $(\cls, U)$ where $\cls$ is a separable unital $C^*$-algebra and  $U$ is a linear map from $\clh$ to $\clh \ot \cls$ such that $\widetilde{U}$ given by $\widetilde{U}( \xi \ot b)=U(\xi) (1 \ot b)$ $(\xi \in \clh$, $b \in \cls$) extends to a unitary element of  $ \clm(\clk(\clh) \ot \cls)$ satisfying the following \\
(i) for every state $\phi$ on $\cls$ we have $U_\phi D=DU_\phi$, wher $U_\phi:=({\rm id} \ot \phi)(\widetilde{U})$;\\
(ii) $({\rm id} \ot \phi) \circ \alpha_U(a) \in ({\cla^\infty})^{\prime \prime}$ $\forall a \in \cla^\infty$ for every state $\phi$ on $\cls$, where $\alpha_U(x):=\widetilde{U}( x \ot 1) {\widetilde{U}}^* $ for $x \in \clb(\clh)$.

In case the $C^*$-algebra $\cls$ has a coproduct $\Delta$ such that $(\cls,\Delta)$ is a compact quantum group and $U$ is a unitary representation  of $(\cls, \Delta)$ on $\clh$, we say that $(\cls, \Delta)$ acts  by orientation-preserving isometries  on the spectral triple.
\edfn

Consider the category ${\bf Q}$ with the object-class consisting of all quantum families of orientation preserving isometries $(\cls, U)$ of the given spectral triple, and the set of morphisms ${\rm Mor}((\cls,U),(\cls^\prime,U^\prime))$ being the set of unital $C^*$-homomorphisms $\Phi : \cls \raro \cls^\prime$ satisfying $({\rm id} \ot \Phi) (U)=U^\prime$. We also consider another category ${\bf Q}^\prime$ whose objects are triplets $(\cls, \Delta, U)$,  where $(\cls,\Delta)$ is a compact quantum group acting by orientation preserving isometries on the given spectral triple, with $U$ being the corresponding unitary representation. The morphisms  are the homomorphisms of compact quantum groups which are also morphisms of the underlying quantum families of orientation preserving isometries. The forgetful functor $F: {\bf Q}^\prime \raro {\bf Q}$ is clearly faithful, and we can view $F({\bf Q}^\prime)$ as a subcategory of ${\bf Q}$.

 Unfortunately, in general ${\bf Q}^\prime$ or ${\bf Q}$ will not have a universal object, as discussed in \cite{qorient}.  We have to reestrict to a subcategory described below to get a universal object in general, though in some cases. Fix a positive, possibly unbounded, operator $R$ on $\clh$ which commutes with $D$ and consider the weakly dense $\ast$-subalgebra $\cle_D$ of $\clb(\clh)$ generated by the rank-one operators of the form $|\xi><\eta|$ where $\xi, \eta$ are eigenvectors of $D$. Define  $ \tau_R(x)= Tr ( R x ),~~x \in \cle_D.$

\bdfn
A  quantum family  of orientation preserving isometries $(\cls, U)$ of $(\cla^\infty, \clh, D)$ is said to preserve the $R$-twisted volume, (simply said to be  volume-preserving if $R$ is understood) if one has  $(\tau_R \otimes  {\rm id} ) (\alpha_U(x))= \tau_R(x)1_\cls$ for all $x \in \cle_D$, where $\cle_D$ and $\tau_R$ are as above.  

If, furthermore, the  $C^*$-algebra $\cls$ has a coproduct $\Delta$ such that $(\cls,\Delta)$ is a CQG and $U$ is a unitary representation  of $(\cls, \Delta)$ on $\clh$, we say that $(\cls, \Delta)$ acts  by ($R$-twisted) volume and orientation-preserving isometries  on the  spectral triple.

 We shall consider the categories ${\bf Q}_R$ and ${\bf Q}^\prime_R$ which are the full subcategories of ${\bf Q}$ and ${\bf Q}^\prime$ respectively, obtained by restricting the object-classes to the volume-preserving quantum families. 
\edfn
The following result is proved in \cite{qorient}.
\bthm
The category ${\bf Q}_R$ of quantum families of volume and orientation preserving isometries has a universal (initial) object, say  $(\widetilde{\clg}, U_0)$. Moreover, $\widetilde{\clg}$ has a coproduct $\Delta_0$ such that $(\widetilde{\clg},\Delta_0)$ is a compact quantum group and $(\widetilde{\clg},\Delta_0,U_0)$ is a universal object in the category ${\bf Q}^\prime_0$.  The representation  $U_0$ is faithful.
\ethm
 
The Woronowicz subalgebra of $\tilde{\clg}$ generated by elements of the form $\{<(\xi \ot 1), \alpha_0(a) (\eta \ot 1)>_{\tilde{\clg}},~a \in \cla^\infty \}$, where $\alpha_0 \equiv \alpha_{U_0}$ and  $< \cdot, \cdot >_{\tilde{\clg}}$ denotes the $\tilde{\clg}$-valued inner product of the Hilbert module $\clh \ot \tilde{\clg}$, is called the quantum group of orientation and volume preserving isometries, and denoted by $QISO^+_R(D)$.

\section{$C^*$-action of $QISO^+_R(D)$}
It is not clear from the definition and construction of $QISO_R^+(D)$ whethe the $C^*$ algebra $\cla$ generated by $\cla^\infty$ is stable under 
$\alpha_0:=\alpha_{U_0}$   in the sense that $({\rm id} \ot \phi) \circ \alpha_0$ maps $\cla$ into $\cla$ for every $\phi$. Moreover, even if $\cla$ is stable, 
the question remains whether $\alpha_0$ is a $C^*$-action of the CQG $QISO^+_R(D)$. 
Although we could not yet decide whether the general answer to the above two questions are afirmative, 
 we have given a number of sufficient conditions for it in \cite{qorient}. In fact, those conditions already cover all classical compact Riemannian manifolds, and 
 many noncommutative ones as well.  In what follows, we shall provide yet another set of sufficient conditions, which will be valid for many interesting spectral triples constructed from Lie group actions on $C^*$ algebras.

 Suppose that there are compact Lie groups $\tilde{G}$, $G$ with a (finite) covering map $\gamma : \tilde{G} \raro G$ (which is group homomorphism), such that the following hold:\\
(a) there is an action $\beta_g$ of $G$ on the von Neumann algera $\cla^{\prime \prime}$ which is strongly continuous w.r.t. the SOT on $\cla^{\prime \prime}$, and moreover, its restriction on $\cla$ is a $C^*$-action, i.e. $g \mapsto \beta_g(a)$ is norm continuous for $a \in \cla$. \\
(b) there exists a strongly continuous unitary representation $   V_{\tilde{g}}$ of $\tilde{G} $ on $\clh$ which commutes with $D$ and $R$, and we also have
   $ V_{\tilde{g}}  a  {V_{\tilde{g}}}^{-1} =  \beta_{g} ( a ) $, where $a \in \cla, \tilde{g} \in   \tilde{G},$ and $ g = \gamma ( \tilde{g} )  $.\\
(c) In the decomposition of the $G$-action $\beta$ on $\cla$ into irreducible subspaces, each irreducible representation of $G$ occurs with at most finite (including zero)  multiplicity.\\ 
  Since $G$ is a Lie group, we choose a basis of its Lie algebra, say $\{ \chi_1, ...,\chi_N\}$. Each $\chi_i$  induce closoble derivations on $\cla$ (w.r.t.  the norm toplogy) as well as on $\cla^{\prime \prime}$ (w.r.t. SOT) which will be denoted by $\delta_i$ and $\tilde{\delta}_i$ respectively. We do have natural Frechet spaces $\cle_1:=\bigcap_{n \geq 1, 1 \leq i_j \leq N}{\rm Dom}(\delta_{i_1} ... \delta_{i_n})$ and $\cle_2:=\bigcap_{n \geq 1, 1 \leq i_j \leq N}{\rm Dom}(\tilde{\delta}_{i_1} ... \tilde{\delta}_{i_n})$, and we assume that \\
(d) $\cle_1$ and $\cle_2$ coincide with $\cla^\infty$.

Let $\hat{G}$ denote the (countable) set of equivalence classes of irreducible representations of $G$, and let $\clv_\pi$ denote the (finite dimensional by assumption) subspace of $\cla^{\prime \prime}$ which is the range of the `spectral projection' $P_\pi$ corersponding to $\pi$, namely $x \mapsto \int_G c_\pi(g) \beta_g(x) dg$, $c_\pi$ being the character of $\pi$. It is easy to see from the assumptions made that elements of $\clv_\pi$ atually belong to $\cla^\infty$, and clearly, this subspace coincides with the range of the restriction of $P_\pi$ on $\cla$. Thus, in particular, the linear span of $\clv_\pi$, $\pi \in \hat{G}$, is norm-dense in $\cla$ as well. 
 
 Now we have the following:
\bthm
Under the above assumptions, ${QISO}^+_R(D)$ has the $C^*$-action of $\clg={QISO}_R^+(D)$ given by the restriction of  $\alpha_0$ on $\cla$ 
\ethm 
{\it Proof:-}\\
 As discussed in the previous section, we consider the reduced CQG $ \rho_h (\clg) $ where $\rho_h$ is the GNS reprsentation of the Haar state, and   denote the Hopf von Neumann algebra (of compact type) ${\rho_h (\clg)}^{ \prime  \prime }$ by $\clm$. 
Clearly, $\alpha_0$ extends to a normal action of $\clm$ on $\cla^{\prime \prime} \subseteq \clb(\clh)$, 
given by $\alpha_0(x)=U_0 (x \ot 1) U_0^* \subseteq \clb(\clh) \ot \clb(L^2(h))$ for $x \in \cla^{\prime \prime}$, 
which decomposes into finite dimensional irreducible subspaces of $\cla^{\prime \prime}$, say $\{ \cla_i, i \in I \}$ ($I$ some index set). Since by assumption, $C(G)$ can be identified with a quantum subgroup of $\clg$, it is clear that each $\cla_i$ is $G$-invariant, i.e. $\beta_g(\cla_i) \subseteq \cla_i$ for all $g$. Thus, we can further decompose $\cla_i$ into $G$-irreducible subspaces, say $\cla_i^\pi$. Clearly, $\cla^\pi_i \subseteq \clv_\pi \subseteq \cla^\infty$. Thus, the restriction of $\alpha_0$ to the linear span of $\clv^{\pi}_i$ 's (say $\clv$), and hence on the $\ast$-algebra generated by $\clv$, is a Hopf algebraic action of $\clm$.  However, by definition of $\clg=QISO^+_R(D)$ it is clear that $\alpha_0|_{\clv^\pi_i}$ is actually a Hopf algebraic action of the CQG $\rho_h(\clg)$.  Moreover, since $\clv^\pi_i$ is finite dimensional, the `matrix coefficients' of $\alpha_0|_{\clv^\pi_i}$ must come from the Hopf algebra $\clg_0$ mentioned in the previous section, on which $\rho_h$ is faithful, and we thus see, by identifying  $\rho_h(\clg_0)$ with $\clg_0$, that $\alpha_0(\clv^\pi_i) \subseteq \clv^\pi_i \ot_{\rm alg} \clg_0$. 

Now it suffices to prove the norm-density of the subspace $\clv$ in $\cla$. To this end, we note that, by the weak density of $\clv$ in $\cla^{\prime \prime}$, the range $\clv_\pi$ of $P_\pi|_{\cla^{\prime \prime}}$ is the weak closure of $P_\pi(\clv)$. But $P_\pi$ being finite dimensional, we have must have that the range coincides with $P_\pi(\clv)$, which is nothing but the linear span of those (finitely many)  $\clv^\pi_i$ which are nonzero , i.e. for which the irreducible of type $\pi$ occurs in the decomposition of $\cla_i$. It follows that $\clv$ contains $\clv_\pi$ for each $\pi$, hence (by the norm-density of ${\rm Span}\{ \clv_\pi,\pi \in \hat{G} \}$ in $\cla$) $\clv$ is norm-dense in $\cla$.
\qed\\
{\bf Examples:}\\
We have following two classes of spectral triples which satisfy our assumptions.\\
(I)
The assumptions of the above theorem are valid (with $R=I$) in case the given spectral triple obtained from an ergodic action of a compact Lie group $G$ on the underlying $C^*$ algebra $\cla$. Let $G$ have a Lie algebra basis $\{ \chi_1, \ldots, \chi_N \}$ and an ergodic $G$-action on $\cla$, it is well-known that $\cla$ has a canonical faithful $G$-invariant trace, say $\tau$, and if we imbed $\cla$ in the corresponding GNS space  $L^2(\tau)$, the operator $i\delta_j$ (where $\delta_j$ is as before) extends to a self adjoint operator on $L^2(\tau)$. Taking $\clh=L^2(\tau) \ot \IC^n$, where $n$ is the smallest positive integer such that the Clifford algebra of dimension $N$ admits a faithful representation as $n \times n$ matrices, we consider $D=\sum_j i\delta_j \ot \gamma_j$, $\gamma_j$ being the Clifford matrices.  The smooth algebra $\cla^\infty$ corresponding to the $G$-action on $\cla$ satisfies our assumption (d) by  Lemma 8.1.20 of \cite{goswami3} (see also \cite{goswami2}), and it has the natural representation $a \mapsto a \ot 1_{\IC^n}$ on $\clh$. If the operator $D$ has a self-adjoint extension with compact resolvent, it is clear that $(\cla^\infty, \clh, D)$ gives a spectral triple satisfying all the assumptions (a)-(d) with $\tilde{G}=G$, $R=I$ and the unitary representation $V$ being $\beta_g \ot I_{\IC^n}$, where $\beta_g$ is the given ergodic action of $G$ on $\cla$, extended naturally as a unitary representation on $L^2(\tau)$. It may be mentioned that the standard spectral triple on the noncommutative tori  arises in this way. \\
(II) Another interesting class of examples satisfying assumptions (a)-(d) come from those classical spectral triples for which the action of the group of orientation preserving Riemannian isometries on $C(M)$ (where $M$ denotes the underlying manifold) is such that any irreducible representation occurs with at most finite multiplicities in its decomposition. It is easy to see that the classical spheres and tori are indeed such manifolds.

\end{document}